\documentstyle{article}

\newtheorem{de}{ \sc Definition}[section]

\newtheorem{pr}[de]{ \sc Proposition}
\newtheorem{th}[de]{ \sc Theorem}
\newtheorem{cor}[de]{ \sc Corollary}
\newtheorem{re}[de]{ \it Remark}
\newtheorem{eq}[de]{\rm}
\newtheorem{qu}[de]{\sc Question}

\def\mpr#1{\;\smash{\mathop{\hbox to 20pt{\rightarrowfill}}\limits^{#1}}\;}
\def\mrp#1{\;\smash{\mathop{\hbox to 20pt{\rightarrowfill}}\limits_{#1}}\;}
\def\mpd#1{\big\downarrow\rlap{$\vcenter{\hbox{$\scriptstyle#1$}}$}}
\def\mdp#1{\llap{$\vcenter{\hbox{$\scriptstyle#1$}}$}\big\downarrow}

\def\lar{\;\smash{\hbox to 16pt{\leftarrowfill}}\;}
\def\ll#1{\mathop{\longleftarrow}\limits^{#1}\limits_{#1}\;} 
\def\lt{\ll{\lar}}
\def\mono{\lhook\joinrel\relbar\joinrel\rightarrow} 
\def\ep{\;\smash{\mathop{\rightarrow\hskip
-7pt\rightarrow}}\;}

\def\qed{\hfill$\Box$\vskip10pt}

\def\IC{I\mkern-1mu C}
\def\IH{I\mkern-1.5mu H}
\def\Tor{\mathop{\rm Tor}\nolimits}
\def\codim{\mathop{\rm codim}}

\def\sqtimes{\,\vbox{\hrule\hbox{\vrule{\sf x}\hskip-0.3pt\vrule}\hrule}\,}

\def\b{\bullet}
\def\Proof{\noindent\hskip 6pt{\it Proof. }}

\font\Bbf=msbm10

\def\Z{{\hbox{\Bbf Z}}}

\def\Q{{\hbox{\Bbf Q}}}
\def\C{{\hbox{\Bbf C}}}
\def\P{{\hbox{\Bbf P}}}
\def\F{{\hbox{\Bbf F}}}
\def\H{{{\cal A}^*}} 

\def\D{{\bf D}_{\rm mix}(B)} 
\def\codim{{\rm codim}\,}

\begin{document}

\author{Matthias Franz\\
{\small\it Section de Math\'ematiques, Universit\'e de Gen\`eve}\\
Andrzej Weber
\footnote{Supported by KBN 1 P03A 005 26 grant, both authors supported 
by EAGER}\\
{\small\it Instytut Matematyki, Uniwersytet Warszawski }}
\date{}
\title{Weights in cohomology and\\ the Eilenberg--Moore spectral sequence}
\maketitle

\begin{abstract} 

We show that in the category of complex algebraic varieties,
the Eilen\-berg--Moore spectral sequence can be endowed
with a weight filtration.
This implies that it degenerates if all spaces involved
have pure cohomology.
As application, we compute the rational cohomology of an algebraic 
$G$-variety $X$
($G$ being a connected algebraic group)
in terms of its equivariant cohomology provided that $H_G^*(X)$ is pure.
This is the case, for example, if $X$ is
smooth and has only finitely many orbits.
We work in the category of mixed sheaves; therefore our results apply
equally to (equivariant) intersection homology.

\end{abstract}

\section{Introduction}
The following spectral sequence was constructed by Eilenberg and Moore
in \cite{EM}:

\begin{th}\label{Eilenberg-Moore}
  Suppose we have a pull-back square of topological spaces
  $$\matrix{A\times_BC&=&D&\rightarrow& A\cr&& \downarrow&&\downarrow\cr
&&C&\rightarrow& B&}$$
with $B$ simply connected and $A\to B$ a fibration.
Then there exists a spectral sequence converging to $H^*(D)$
with $$E^{-p,q}_2=
\Tor^{q,H^*(B)}_p(H^*(A),H^*(C)).$$
\end{th}

Cohomology is taken with rational coefficients.
The torsion product is 
a functor of homological type and therefore it is 
denoted by $\Tor_p^{H^*(BG)}(-,-)$.
It is positively graded. Additionally
it has an internal grading.
Its degree--$q$ piece is denoted by
$\Tor^{q,H^*(BG)}_p(-,-)$.
According to L.~Smith~\cite{Sm1, Sm2}, the
entries
of the spectral sequence are the cohomology groups of 
certain spaces and the differentials are induced by maps between
them. 
From this Smith deduced that the
sequence inherits all the structure of cohomology, such as an action
of the Steenrod algebra. We are interested in the weight filtration
of the rational cohomology of algebraic varieties. Since we want to prove some
results for singular spaces, we deal with
intersection cohomology as well.                                  
We
show
that Smith's construction can be carried out in the category of
algebraic varieties, or rather
in the derived category
of sheaves on algebraic varieties. Therefore
the Eilenberg--Moore spectral sequence can be endowed with a weight
filtration which extends the filtration constructed by Deligne \cite{D2}. 
The Eilenberg-Moore spectral sequence can
equally
be obtained from a
filtration of the bar complex $B(\Omega^*(A),\Omega^*(B),
\Omega^*(C))$,
which is quasi-isomorphic to $\Omega^*(A\times_B C)$
(where $\Omega^*(-)$ is a complex of differential forms
on a hyperresolution computing $H^*(-)$, see \cite{D2}).
According to Deligne, these complexes live in the category
of mixed Hodge complexes. By \cite[\S3]{Hain} and
\cite[\S3]{Hain-Zucker},
the bar complex inherits the mixed Hodge structure, hence also
the entries of the Eilenberg--Moore spectral sequence.
We only study the weight filtration, but we are
interested
in intersection cohomology as well.
We prove our results in the category of mixed sheaves of Saito
(\cite{Sa1, Sa2}) or Beilinson-Bernstein-Deligne (\cite{BBD}). 
We formulate our first
theorem in the following way:

\begin{th} Let $B$ be a simply connected 
complex algebraic variety.  Let $F$
and $G$ be bounded below mixed sheaves over~$B$. 
Suppose that $F$ has constant cohomology sheaves. Then there is a
spectral sequence with
$$E_2^{-p,q}=\Tor^{q,H^*(B)}_p\bigr(H^*(B;F),H^*(B;G)\bigl)$$
and converging to $H^*(B;F\otimes G)$. 
The entries of the spectral sequence are endowed with weight
filtrations. The differentials preserve them.
\end{th}

Here $H^*(B)$ denotes the cohomology of $B$ with coefficients in the
constant sheaf $\Q$.
The choice of the category
of mixed sheaves
is motivated by the following:
In contrast to the topological situation, we would have to deal with
simplicial varieties instead of just varieties. That is because
in the construction of Smith's resolution on the geometric level a
quotient of varieties appears. Such a quotient is no longer an algebraic
variety. In the world of simplicial varieties we can replace
the quotient by a cone construction, and every step of
Smith's original proof can be imitated.
But we are concerned with intersection homology as well,
and we would
therefore
have to
introduce a simplicial intersection sheaf. We find it
rather unnecessary to develop a theory of intersection cohomology
for simplicial varieties. 
On the other hand it is much
easier and more general to prove the result for sheaves. If one looks
carefully, one sees that Smith's construction is carried out in the
stable category of topological spaces over $B$ (at the end of his
argument he desuspends the spectral sequence). The rational stable category of
topological spaces over $B$ is nothing but the derived category of
sheaves over $B$. The sheaves coming from algebraic geometry carry an
additional structure: the weight filtration. In a purely formal way
we deduce results from the fact that all maps strictly preserve
the filtrations. In many
situations the higher differentials of the Eilenberg-Moore spectral
sequences vanish because otherwise they would mix weights.
Hence, $E_2=E_\infty$ in these cases.
We note that if a variety is smooth, one could work
with Hodge complexes as in \cite{D2} instead of mixed sheaves.

Using the same notation as before, our main result reads as follows:

\begin{th} \label{D} If $H^*(B)$, $H^*(B;F)$ and
$H^*(B;G)$ are pure
(i.e., the degree-$n$ cohomology group is entirely of weight $n$),
then $E^{-p,*}_r$ is pure for all $r$~and~$p$,
and the Eilenberg--Moore spectral sequence degenerates
on the $E_2$~level.
The resulting filtration is the weight filtration
$$W_\nu H^{n}(F\otimes G)=\bigoplus_{q-p=n,\;q\leq
\nu}\Tor^{q,H^*(B)}_p(H^*(B;F),H^*(B;G))\,.$$ 
\end{th}

Suppose we have a pull-back square
\begin{eq}\label{diag} \hfil $\matrix{A\times_BC&=&D&\mpr{\tilde g}&A\cr
&&\mdp{\tilde f}&&\mpd{f}\cr
&&C&\mrp{g}&B\cr}$\end{eq}
\noindent of algebraic varieties, where $B$ is smooth and simply
connected and the 
map $A\rightarrow B$ 
is a fibration.
(It is enough to assume that the map $f$ is a topological
locally trivial fibration.)
When we apply Theorem \ref{D} to the push forwards
$F=Rf_*\IC_A$~and~$G=Rg_*\IC_C$ of
the intersection sheaves, we obtain:

\begin{th}\label{C} If $H^*(B)$, $\IH^*(A)$  and $\IH^*(C)$ are pure, then
$$\IH^n(D)=\bigoplus_{q-p=n} \Tor^{q,H^*(B)}_p(\IH^*(A),\IH^*(C))\,.$$
The sum of terms with $q\leq \nu$ coincides with $W_\nu \IH^*(D)$.\end{th}

Let $G$ be a linear algebraic group. We assume that $G$ is connected.
Our goal is to study the weight structure in the cohomology and 
equivariant cohomology of algebraic $G$-varieties. 
If the equivariant cohomology is pure,
it determines the non-equivariant cohomology additively:
  
\begin{th}\label{B} If 
the rational equivariant cohomology 
$H^*_G(X)=H^*(EG\times_GX)$
is pure, then
the rational cohomology of $X$ is given additively by:
$$H^n (X)=\bigoplus_{q-p=n} \Tor^{q,H^*(BG)}_p(H^*_G(X),\Q)\,.$$
The sum of terms with $q\leq \nu$ coincides with $W_\nu 
H^*(X)$.\end{th}

Theorem \ref{B} follows from 
Theorem \ref{C} by approximating $BG$ in the pull-back square
$$\matrix{ X&\rightarrow& EG\times_GX\cr
\downarrow&&\downarrow\cr
point&\rightarrow&BG&\cr}\,.$$ 
In the special case where $X=G/H$ is a homogeneous spaces we recover a
result of Borel \cite[Th. 25.1 or 25.2]{Bo}

An analogous theorem can be formulated for intersection cohomology
or Borel--Moore homology.
It remains to say when equivariant cohomology is pure.
Without difficulty, we find (see Propositions \ref{fi} and \ref{add}):

\begin{th} If a $G$-variety $X$ is smooth 
and has only finitely many orbits,
then $H^*_G(X)$  
is pure. 
The rational cohomology of $X$ is given additively by:
$$H^*_G(X)=\bigoplus H^{*-2c}(BH)\,,$$
where the sum is taken over all orbits ${\cal O}=G/H\subset X$,
and $c=\codim\,\cal O$.
\end{th}  

For singular varieties the result holds for equivariant Borel--Moore
homology as 
defined in \cite[Section 2.8]{EG}.
The multiplicative structure of $H^*_G(X)$ is harder. It involves
Chern classes of normal bundles of orbits. Except for special cases
like toric varieties or some other
spherical varieties, we do not have a satisfactory description.

In order to apply Theorem \ref{B} to intersection cohomology of singular
varieties, we have to show that 
$\IH^*_G(X)$ is pure. An important class of varieties having this property
is that of spherical varieties.
One can prove it by using
a local description of singularities as in~\cite{BJ2}.

{\it
Standing assumptions:}
Unless stated otherwise,
all cohomology  groups   are   taken   with 
rational coefficients.
All algebraic varieties are defined over the complex
numbers.
We consider only algebraic actions of linear algebraic groups.
\medskip

We would like to thank the referee for his great help and care.

\section{Preliminaries}

\subsection{Weight filtration}
The weight filtration in the cohomology of algebraic varieties
can be constructed in various ways: either using arithmetic methods,
\cite[\S6]{D1} 
or through analytic methods, \cite{D2}. To tackle 
not 
necessarily constant coefficient systems
one should work with the mixed Hodge modules of
M.~Saito, \cite{Sa1, Sa2}.  
Our results will follow formally from the existence of a weight filtration.
Therefore, instead of going into various constructions, we will list its
properties.  Let  $X$  be  a  smooth  variety.  We  consider
cohomology with coefficients in $\Q$.
The weight filtration
$$0= W_{k-1}H^kX\subset W_kH^kX\subset W_{k+1}H^kX\subset \dots 
\subset W_{2k}H^kX=H^kX$$ satisfies the following conditions:
\begin{description}
\item[1.] $W_kH^kX=H^kX$ if $X$ is complete.
(We say that the cohomology is \emph{pure}.)

\item[2.] Let $f\colon X\rightarrow Y$ be an algebraic map. The induced map
in cohomology strictly preserves weight, i.e.,
$f^*W_\nu H^kY=W_\nu H^kX\cap f^*H^kY$.

\item[3.] The weight filtration is strictly preserved by the maps in the
Gysin (localization) sequence
$$\cdots\mpr{} H^kX \mpr{j^*} H^kU \mpr{\delta} H^{k+1-2c}Y\mpr{i_!}
H^{k+1}X\mpr{}\,\cdots. $$
Here $i\colon Y\mono X$ is a smooth subvariety of codimension $2c$, and
$j$ is the inclusion~$U=X\setminus Y\mono X$. The boundary map
 $\delta$ lowers weight by~$2c$ and $i_!$ raises it by~$2c$.

\item[4.] The weight filtration is defined for relative cohomology.
The long exact sequence of a pair strictly preserves it.

\item[5.] The weight filtration is also defined for simplicial varieties.
This time weights smaller than $k$ can appear in $H^k(X_\b)$.
By hyperresolution of singularities, the cohomology of singular
varieties can be endowed with a weight structure.
\end{description}

We deal with singular varieties, therefore we have to consider 
not only constant sheaves, but also the ``mixed sheaves'' of~
\cite[\S5, p.~126]{BBD}. In
particular, we compute cohomology with coefficients in the intersection
sheaf $\IC_X$, that is, intersection cohomology~\cite{GM}. 
In the main part of~\cite{BBD} the varieties are actually defined over
a field of finite characteristic. 
For complex varieties we are allowed to work with ``sheaves of
geometric origin'' 
(\cite[\S6.2.4]{BBD}).
For such sheaves it is possible to find a good reduction and
apply the results which are valid in finite characteristic. 
In particular we have $\IH^*(X_\C);\Q_\ell)\simeq
\IH^*(X_{\overline{\F}_q};\Q_\ell)$. 
(We recall that for varieties over a finite field we have
to use coefficients in $\Q_\ell$. The formalism of \cite{D1} which
allows to compute the cohomology with coefficients in a complex of
sheaves generalizes the construction of \'etale cohomology.)
This way we transport the weight filtration to the intersection cohomology of
complex varieties. 
The theory of M.~Saito also applies: an
intersection sheaf is an object 
in the category of mixed Hodge modules. 
The weight filtration in $\IH^k(X)$ also starts with 
$W_k\IH^kX$ as in the smooth case. 
It turns out that the weight filtration is defined over rational numbers.
We will use again just formal properties of mixed sheaves. 
A suitable category to work with is the
category of mixed Hodge modules over a base $B$, i.e. $\D$.
The
properties 3~and 4~above can be extended to the following:
\begin{description}
\item[6.]  For a distinguished triangle of mixed sheaves 
$$A\mpr{}B\mpr{}C\mpr{}A[1]$$ the maps in the long
exact sequence 
$$\cdots\to H^k(X;A) \to H^k(X;B) \to H^k(X;C)\to
H^{k+1}(X;A)\to\,\cdots$$
strictly preserve weights. 
\end{description}

\subsection{Varieties associated with a group action}

Let $G$ be a connected linear algebraic 
group. The classifying space of~$G$ as a
simplicial variety has already been considered by Deligne~\cite{D2}:
$$BG_\b=\left\{pt\; ^{\lar}_{\lar} G \lt G\times G\; \cdots\right\}$$
(the arrows are multiplications or forgetting the edge factors).
In an appropriate category, $BG_\b$ represents the functor which
associates to~$X$ the isomorphism classes of Zariski-locally trivial
$G$-bundles over~$X$.  

Another algebraic model
$BG_{\acute et}$ was described
by Totaro~\cite{To}, see also~\cite[\S4.2]{MV}. It
classifies $G$-bundles which are \'etale-locally trivial.
Totaro's~$BG_{\acute et}$ has the same cohomology as the simplicial~$BG_\b$. 
The spaces $EG_{\acute et}$~and~$BG_{\acute et}$ are infinite-dimensional;
one therefore has to approximate them by $U$ and $U/G$
where $U\subset V$ runs over all open
$G$-invariant subsets in representations 
$G\to GL(V)$, such that the geometric quotient $U/G$ exists and $U\to U/G$
is a $G$-bundle (\'etale locally trivial).
We prefer to work with this model. 
Both models have isomorphic
cohomology.

The Borel construction of~$X$
is the simplicial variety
$$[X/G]_\b=(EG\times_G X)_\b=\left\{X\; ^{\lar}_{\lar} X\times G \lt X\times
G\times G\; \cdots\right\}$$ 
(the arrows are the action, multiplications or forgetting the last factor).
Of course, $[pt/G]_\b=BG_\b$.

Note that in Totaro's model
the space~$U\times_G X$ might not be an algebraic 
variety.
But it is so 
if all
orbits admit  $G$-invariant quasi-projective neighbourhoods.
See the discussion in~\cite{We}.

Equivariant cohomology $H^*_G(X)$ is defined as the cohomology of 
$EG\times_GX$. It does not matter which model of $EG\to BG$ we use. 
For a fixed degree~$i$ we have
$H^i_G(X)=H^i(U\times_GX)$ if 
the codimension of $V\setminus U$ is sufficiently large.

\subsection{The weight filtration in~$H^*(BG)$ and in~$H^*(G)$}

The following observations were made by Deligne~\cite{D2}, \S9. 
  
\begin{th}\label{BG-pure} The cohomology of~$BG$ is pure, i.e.,
$$W_{k-1}H^k(BG)=0\,,\qquad W_kH^k(BG)=H^k(BG)\,.$$
\end{th}
  
For example, $H^*(B\C^*)\simeq H^*(\P^\infty)$ is a polynomial algebra on
one pure generator of degree~$2$.
  
Now let $P^\b=P^1\oplus P^3\oplus P^5\oplus\dots$ be the space of
primitive elements of the Hopf algebra
$H^*G$, so that
$H^*G=\bigwedge P^\b$.
  
\begin{th} Let $k>0$. Then $W_kH^kG=0$ and $W_{k+1}H^kG=P^k$ (which is 0
if $k$ is even).\end{th}

Hence one can filter the cohomology of~$G$ by ``complexity'':
$C_a H^*(G)=\bigwedge\nolimits^{\leq a}P^\b$.
Then $C_a H^*(G)\cap H^k(G)=W_{a+k}H^k(G)$.

\section{Filtration in the Eilenberg--Moore spectral sequence}
Let us come back to the pull-back diagram \ref{diag} and the
associated  Eilenberg-Moore spectral sequence.
We want to endow its entries 
with a weight filtration. The goal of this section is to repeat
Smith's construction (as presented in \cite{Sm2}).
We find it convenient to consider the category
of sheaves over~$B$ instead of the category of spaces over~$B$
as in Smith's papers. More precisely, we work in the category~$\D$
of mixed sheaves
of~\cite{Sa1,Sa2} or \cite{BBD}. Instead of a map $f\colon A\rightarrow B$
we consider the sheaf~$Rf_*\IC_A$ since~$\IH^*(A)=H^*(B;Rf_*\IC_A)$.
The condition: {\it $f\colon A\to B$ is a fibration with the fibre
$X$} is
replaced by:
{\it $Rf_*\IC_A$ has constant cohomology sheaves
with stalks $\IH^*(X)$}. Since we assume $B$ to be simply connected
there is no need to consider locally constant cohomology sheaves.
It is clear that if $f\colon A\to B$ is a fibration which is locally
trivial with respect to the classical topology, then $Rf_*\IC_A$ has
constant cohomology sheaves. Note that in this case the weight
filtration is constant, \cite[Proposition 6.2.3]{BBD}, although the Hodge
structure might vary.
Suppose that $B$ is smooth. 
We claim that the intersection cohomology of the pull-back is the
cohomology of~$B$
with coefficients in the tensor product $Rf_*\IC_A\otimes Rg_*\IC_C$:
$$\IH^*(D)\simeq H^*(B;Rf_*\IC_A\otimes Rg_*\IC_C)$$
or, more precisely, that
$\bar f^*\IC_C\otimes \bar g^* \IC_A$ is quasi-isomorphic to $\IC_D$.
Indeed, 
\begin{itemize}
\item this sheaf is constant when restricted to the regular part of
$D$, which is equal to $C_{\rm reg}\times_BA_{\rm reg}$;
\item locally, 
for an open set $U\subset C$ (in the classical topology) 
over which the fibration is
trivial, i.e., $\bar f^{-1}(U)=U\times X$, we have $(\bar f^*\IC_C\otimes
\bar g^* \IC_A)\simeq \IC_U\otimes \IC_X$.
\end{itemize}

Therefore our goal is to construct a spectral sequence converging to
cohomology with coefficients in a tensor product of mixed sheaves.
The same applies to constant sheaves instead of intersection sheaves.
For~$B=BG$ one can generalize our construction to equivariant
sheaves in the sense of Bernstein-Lunts, \cite{BL}. 

We set $\H=H^*(B)$, 
and we write $H^*F$ instead of~$H^*(B;F)$  for a
mixed sheaf~$F$.

\begin{th} Let $B$ be a simply connected 
 complex
algebraic variety.  Let $F$
and $G$ be bounded below mixed sheaves over~$B$. 
Suppose that $F$ has constant cohomology sheaves. Then there is a
spectral sequence with
$$E_2^{-p,q}=\Tor^{q,\H}_p(H^*F,H^*G)$$
and converging to $H^*(F\otimes G)$. 
The entries and differentials of the spectral sequence lie in the category
${\bf D}_{\rm mix}(point)$.
\end{th}

\Proof Let $\pi\colon B\times B\to B$
be the projection onto the first factor and
let $\Delta\colon B\to B\times B$ be the diagonal. For a mixed 
sheaf
$F$ over $B$ 
we define $QF=R\pi_*\pi^*F$. It comes with a map $QF\to F$
constructed in the following way:
\begin{enumerate}
\item $F\to F=R\pi_*R\Delta_*F$ is the identity map over $B$,
\item $\pi^*F\to R\Delta_*F$ is the adjoint map,
\item applying $R\pi_*$, we obtain $R\pi_*\pi^*F\to
R\pi_*R\Delta_*F=F$. 
\end{enumerate}
The sheaf $QF$ is isomorphic to $R\epsilon_*\Q_B$\sqtimes$F$,
where $\epsilon\colon B\to point$. 
The cohomology of the stalk ${\cal H}^*_x(QF)$ is equal to
$\H\otimes {\cal H}^*_x(F)$.
We note that $H^*QF$ is a free
$\H$-module and that the induced map $H^*QF\rightarrow H^*F$ is surjective.
Indeed, $H^*QF=\H\otimes H^*F$ and the map to $H^*F$ is the action of
$\H$ on $H^*F$. Let us
summarize the ingredients which L.~Smith has used to construct the
spectral sequence (most of them are listed in~\cite{Sm1}):

\begin{enumerate}
\item \label{property1}
  For all $F\in Ob(\D)$ the map $QF\rightarrow F$ satisfies the
  properties 
\begin{description}
\item{a)} $H^*QF\rightarrow H^*F$ is surjective,

\item{b)} $H^*QF$ is a free $\H$-module,

\item{c)} if $H^qF=0$ for $q<a$, then $H^qF\rightarrow H^qQF$ is an
isomorphism for $q<a+2$. (This condition is necessary for
convergence, see~\cite{Sm2}, p.~42.)
\end{description}

\item If $F$ is a fibration (i.e., has constant cohomology sheaves), 
then so is $QF$.

\item If $H^*F$ is of finite type, then $H^*QF$ is of finite type as well.

\item \label{property4}
If $F$ is a fibration and $H^*F$ a free $\H$-module, then the natural map 
$$H^*F\otimes_\H H^*G\rightarrow H^*(F\otimes G)$$ 
is an isomorphism for all~$G\in Ob(\D)$.

\item \label{property5}
Let $F$,~$G\in Ob(\D)$ with $F$~being a fibration.
If $H^qF=0$ for~$q<a$ and $H^qG=0$ for~$q<b$,
then $H^q(F\otimes G)=0$ for~$q<a+b$.
(This is again necessary for convergence, compare~\cite{Sm2},
p.~43 and Proposition~A.5.2.)
\end{enumerate}

Properties (\ref{property4})~and~(\ref{property5})
follow from the Grothendieck spectral sequence
$$E^{p,q}_2=H^p({\cal H}^qF\otimes G)\Rightarrow H^{p+q}(F\otimes
G)\,.$$ 

Now we mimic Smith's construction of a free resolution of
$H^*F$. We inductively define
\begin{itemize}
\item $F_0=F$;
\item $Q_p=QF_p$;
\item $F_{p+1}$ is the fibre of $Q_p=QF_p\rightarrow F_p$, i.e., it fits
into the distinguished triangle 
$$\matrix{Q_p&
\hskip-15pt\;\smash{\mathop{\hbox to 30pt{\rightarrowfill}}}\;\hskip-15pt 
& F_p.\cr
&\nwarrow\phantom{^\star}\swarrow&_{[+1]}\hfill\cr
&F_{p+1}&\cr}$$
\end{itemize}
We thus obtain a free resolution of~$H^*F$ coming on the level of sheaves
from
$$
\def\ma{\hskip-15pt\;\smash{\mathop{\hbox to
30pt{\leftarrowfill}}}\;\hskip-15pt } 
\def\mb{\;\smash{\mathop{\hbox to
30pt{\rightarrowfill}}\limits^{[+1]}}\;} 
\def\st{{^\star}}
\matrix{
F=F_0&\mb&F_1&\mb&F_2&\dots& \mb&F_p&\mb&F_{p+1}\cr
&\nwarrow\st\swarrow&_\bigtriangleup&\nwarrow\st\swarrow& &
_\bigtriangleup&\nwarrow\st\swarrow&_\bigtriangleup&\nwarrow\st
\swarrow&\cr 
& Q_0&\ma&Q_1&\dots\;&\ma&Q_p&\ma&Q_{p+1}&\dots\cr
}
$$
Here $\star$ denotes a distinguished triangle
and $\bigtriangleup$ a commutative triangle.
The sheaves~$Q_p$ have constant cohomology.
We tensor the entire diagram by~$G$. The triangles remain distinguished. Now
we apply $H^*(-)$.
The Eilenberg--Moore spectral sequence is obtained from the exact couple
$$\matrix{\bigoplus_pH^*(F_p\otimes G)&
\mathop{\hbox to 50pt{\rightarrowfill}}\limits^{[+1]}
& \bigoplus_pH^*(F_p\otimes G).\cr
\hfill\nwarrow& &\swarrow\hfill\cr
&\bigoplus_pH^*(Q_p\otimes G)&\cr}$$
This is an exact couple in the abelian category of vector spaces with
an action of the Frobenius automorphism. The maps strictly preserve the
weight filtration. If we worked with Saito's category of mixed Hodge
sheaves, the cohomology groups would lie in the abelian category of
mixed Hodge structures; here the maps also strictly preserve the
weight filtration.                      
The elements of the couple are bigraded by~$p$
and by the internal degree~$q$. The map~$\rightarrow$ is of
bidegree~$(1,1)$, the map~$\swarrow$ of bidegree~$(0,0)$
and the map~$\nwarrow$ of bidegree~$(1,0)$.
The resulting spectral sequence has
$E^{-p,q}_1=H^q(Q_p\otimes G)$.
By property~(\ref{property4}) we have 
$$E^{-p,*}_1=H^*(Q_p\otimes G)=H^*Q_p\otimes_\H H^*G\,.$$
The $E_2$~term 
is
$$E_2^{-p,*}=H_p(H^*(Q_\b\otimes G))=\Tor_p^\H(H^*F,H^*G)\,.$$

\begin{re}\rm The properties of this spectral sequence were studied in detail 
by
L.~Smith in
\cite{Sm2}. His arguments fit perfectly into the formalism of the stable
homotopy category over $B$. In the final step he desuspends
his spectral sequence
\cite[p.~38]{Sm2}. 
The formalism of the derived category of sheaves
over $B$ is even better for this purpose.\end{re}

The convergence of the Eilenberg--Moore spectral sequence is
guaranteed by the following.
Suppose that the mixed sheaves $F$ and $G$ have cohomology concentrated
in non-negative degrees.
Then the resulting spectral sequence lies in the second quadrant.
Inductively one checks that $H^qF_p=0$ for $q<2p$. Indeed, by 
property (\ref{property1}c) the map
$H^qQ_p\to H^qF_p$ is an isomorphism for~$q<2p+2$ and a 
surjection
for~$q=2p+2$. Therefore $H^qF_{p+1}=0$ for $q<2p+2$.
Now by property~(\ref{property5}) $H^q(Q_p\otimes G)=0$ for~$q<2p$ 
(since $H^q(Q_p)=H^q(F_p)=0$). 
We find that
the terms $E^{-p,q}_2$ may be nonzero only for~$0\geq q\geq 2p$.
\qed

The entries of the spectral sequence are cohomologies of objects 
in ${\bf D}_{\rm mix}(point)$, therefore:

\begin{cor} The Eilenberg--Moore spectral sequence inherits a
weight filtration.\end{cor}

Let us assume that $\H$, $H^*F$ and $H^*G$ are pure.
Then $H^*Q_0=\H\otimes H^*F$ is pure.
The cohomology~$H^*F_1$ is the kernel of~$H^*Q_0\ep H^*F$,
hence also pure. Reasoning inductively, we find that 
$H^*Q_p$ is pure for all~$p$.
Hence $E^{-p,q}_1=(H^*Q_p\otimes_\H H^*G)^q$ 
is pure of weight~$q$.
Therefore the differentials vanish from $d_2$ on.
We obtain:

\begin{th} \label{pur} If $\H$, $H^*F$ and
$H^*G$ are pure, then $E^{-p,*}_r$ is pure. Consequently, the
Eilenberg--Moore spectral sequence degenerates and
$$W_\nu H^{n}(F\otimes G)=\bigoplus_{q-p=n,\;q\leq
\nu}\Tor^{q,\H}_p(H^*F,H^*G)\,.$$ 
\end{th}

\begin{re}\rm Our construction is explicit.
The spectral sequence can be deduced from the simplicial sheaf~%
$F\otimes Q_\bullet G$. The term
$E_1^{-p,*}=H^*(Q_p\otimes G)$ can be
identified with the bar complex 
$${\mathop{\rm Bar}}_p(H^*F,\H,H^*G)=H^*F\otimes
(\overline {\cal A}^*)^{\otimes p}\otimes H^*G\,.$$
Here $\overline {\cal A}^*=\overline H^*(B)$ is the reduced cohomology of the
base.   
\end{re}

\begin{re} \rm We will apply the Eilenberg--Moore spectral sequence
to sheaves over $BG$. We can extend Theorem~\ref{pur} to this case
by approximating $BG$ the way Totaro does.\end{re}

\section{Applications} 

\subsection{
 Rational cohomology of principal bundles}

Let $G$ be a connected linear algebraic group
and let $X$ be a complete  complex
algebraic variety.
Let $P\to X$ be a principal algebraic $G$-bundle
which is \'etale locally trivial.
It does not have to be Zariski-locally trivial, but 
the map $P\to X$ is affine.
In general a classifying map $X\to BG$ might not exist in the
category of algebraic varieties.
Instead, we use a construction of
\cite[Proof of Theorem 1.3]{To}:
Let $V$ be a representation of $G$ and $X_V=P\times_GV\to X$ the
associated vector bundle. 
The algebraic variety $X_V$ is homotopy 
equivalent to $X$. Let $U$ be an open $G$-invariant subset of $V$ with 
a quotient $U/G$ which approximates $BG_{\acute et}$.
The intersection cohomology of the open set $X_U=P\times_GU\subset X_V$
approximates $\IH^*(X)$ as $\codim 
(V\setminus U)\to \infty$. 
Similarly  $\IH^*(P\times U)$ approximates $\IH^*(P)$. 
Now there is a pull-back
  $$\matrix{P\times U&\rightarrow& BE_{\acute et}
\cr \downarrow&&\downarrow\cr
X_U&\rightarrow& BG_{\acute et}&}$$
Set $B=BG_{\acute et}$, $A=EG_{\acute et}$ and $C=\lim X''$. The Eilenberg--Moore
spectral sequence allows to compute the intersection cohomology of
$P$. 
Since $\IH^*(X)$ is pure (as is $H^*(EG)$, of course),
we obtain:

\begin{cor} 
$$W_\nu \IH^{n}(P)=\bigoplus_{q-p=n,\;q\leq
\nu}\Tor^{q,H^*(BG)}_p(\IH^*(X),\Q)\,.$$ 
\end{cor}

The approximation is justified by the following proposition,
which follows simply from the fact that Tor can be computed from
the bar-resolution.

\begin{pr}
Let $A^*$ be an inverse limit of graded rings $A^*_i$ ($i\in \Z$) and let
 $M^*$ and $N^*$ be
graded $A^*$-modules which are inverse limits of $A^*_i$-modules
$M^*_i$ and $N^*_i$.
We assume that the grading is nonnegative and
that for each degree~$n$, the limits stabilize
(i.e., $A^n=A^n_i$ for large~$i$ and analogously for $M^*$ and $N^*$).
Then for a fixed $q$ and sufficiently large $i$
$$
  {\rm Tor}_*^{q,A^*}(M^*,N^*)={\rm Tor}_*^{q,A^*_i}(M^*_i,N^*_i).
$$
\end{pr}

\subsection{From equivariant to non-equivariant cohomology}

Let $G$ be as before and $X$ a $G$-variety.
Set $B=BG$, $A=EG\times_GX$ and~$C=point$.
We list examples
to which Theorem \ref{B} applies, that is,
for which $H^*(EG\times_GX)=H^*_G(X)$ or 
$\IH^*(EG\times_GX)=\IH^*_G(X)$ are pure.
Let us start with some simple observations. 

\begin{pr}\label{fi}
Let $X\to Y$ be a fibration (in the classical topology) 
with fibre, a homogeneous $G$-space.
If $Y$ is complete, then $\IH^*_G(X)$ is pure.
\end{pr}

\Proof Consider the fibration~$EG\times_GX\rightarrow Y$. The fibres
are homotopy equivalent to~$BH$, where $H$ is the stabilizer of a point in~%
$X$. Then by the Serre spectral sequence we find that
$\IH^*(EG\times_GX)$ is pure.
(Note that this spectral sequence degenerates.)\qed

\begin{pr}\label{fi2}
Suppose $X$ is equivariantly fibred over a homogeneous $G$-space.
If the fibres are complete, then $\IH^*_G(X)$ is pure.
\end{pr}

\Proof
The fibration is of the form~$X\to G/H$. Let $F$ be a fibre, which is complete.
Consider the fibration~$EG\times_GX\rightarrow
EG\times_GG/H\simeq BH$. The fibres are again isomorphic to~$F$.
The Serre spectral sequence has $E_2^{p,q}=H^p(BH;\IH^q(F))$. 
(The coefficients might be twisted.)
We have
$$H^p(BH;\IH^q(F))=(H^p(\widetilde{BH})\otimes \IH^q(F))^{\pi_1(BH)}
=(H^p(BH^0)\otimes (\IH^q(F)))^{H/H^0}\,,$$
where $H^0$ is the identity component of~$H$.
We see that $E^{p,q}_2$ is pure of weight $p+q$. Therefore
$\IH^*_G(X)$ is pure. Moreover, the spectral sequence degenerates. 
\qed

Now suppose that $X$ is smooth. We show that purity is additive in this case.

\begin{pr}\label{add} Suppose that $X$ is smooth and has a $G$-equivariant
stratification~$\{S_\alpha\}$ such that each~$H^*_G(S_\alpha)$ is
pure. Then
\begin{equation}
  \label{cohom-strat}
  H^*_G(X)= \bigoplus_p\bigoplus_{\codim S_\alpha=p}
              H^{*-2p}_G(S_\alpha)(p) 
\end{equation}
is pure.
\end{pr}

Here $(p)$ denotes the shift of weights by~$2p$.

\Proof Consider the filtration of $X$ by codimension of strata:
$$U_p=\bigcup_{\codim S_\alpha\leq p}S_\alpha\,.$$
The resulting spectral sequence for equivariant cohomology has
$$E^{p,q}_1=H^{p+q}_G(U_p,U_{p-1})\,.$$
By the Thom isomorphism,
$$E^{p,q}_1=\bigoplus_{\codim S_\alpha= p} H^{q-p}_G(S_\alpha)\,.$$
The Thom isomorphism lowers weight by~$2p$. Hence,
$E^{p,q}_1$ is pure of weight~$p+q$
because $H^{q-p}_G(S_\alpha)$ is pure of weight~$q-p$. We conclude that the
spectral sequence degenerates,
which was to be shown.
\qed

By Theorem~\ref{BG-pure},
Proposition~\ref{add} applies in particular to varieties with finitely
many orbits. Since the right hand side of~(\ref{cohom-strat})
is the $E_1$~term of a cohomology
spectral sequence, it carries a canonical product.
For a toric variety~$X$
defined by a fan~$\Sigma$, this ``is'' the product
in~$H_G^*(X)$ (which equals the Stanley--Reisner ring of~$\Sigma$).
The same holds for complete symmetric varieties~\cite[Thm.~36]{BDCP} and,
more generally, for toroidal varieties~\cite[\S 2.4]{Br}. 

\medskip

Now we switch to the singular case.

\begin{re}\rm
First we note that
our
results can be generalized straightforwardly if
we replace the constant sheaf by the dualizing sheaf. This way
one compares Borel--Moore homology with equivariant Borel--Moore
homology, defined in \cite{EG}. We will not develop this remark
here. \end{re}

Let us make some comments about intersection sheaves,
which compute intersection cohomology \cite{GM}.
If the local intersection sheaf is pointwise pure,
we can argue as in the proof of Proposition~\ref{add}. 
This condition is satisfied if
the singularities are quasi-homogeneous (as 
in~\cite{DL}, for instance). Other examples are obtained
from the decomposition theorem (\cite[Th\'eorem\`e 6.2.5]{BBD}, \cite{Sa3}):

\begin{pr} If $X$ admits an equivariant resolution~$\widetilde X$ with pure
equivariant cohomology, then the equivariant intersection cohomology of~$X$
is pure.\end{pr}

\Proof Indeed, $\IH_G^*(X)$ injects into $H_G^*(\widetilde X)$
by the decomposition theorem,  and the inclusion preserves weights.\qed

There are lots of examples of $G$-varieties admitting a resolution which
can be 
stratified by~$S_\alpha$ as in~\ref{add}. 
If $X$ is a spherical variety, one can prove purity without appealing
to the decomposition theorem. Here a precise analysis of spherical
singularities (\cite{BJ2} \S3.2) shows that the intersection
sheaf is pointwise pure. The stalk ${\cal H}^*_x\IC_X$ coincides with the
primitive part of intersection cohomology of some projective spherical
variety of lower dimension. 
We refer the reader to~\cite{BJ1}, Theorem 2 and~\cite{BJ2},
Theorem 3.2.
Filtering $X$ by codimension of orbits, we find that
$$\IH^*_G(X)=\bigoplus_\alpha H^*(BG_\alpha;{\cal H}^*_{x_\alpha}\IC_X)\,.$$
The sum is taken over all orbits $G\cdot x_\alpha=G/G_\alpha$.
The decomposition follows from the purity of the summands.
The coefficients may be twisted.
The intersection cohomology localized at the orbit $G\cdot x_\alpha$ is 
equal to
$$H^*(BG_\alpha;{\cal H}^*_{x_\alpha}\IC_X)=
H^*(BG^0_\alpha;{\cal H}^*_{x_\alpha}\IC_X)^{G_\alpha/G_\alpha^0}\,.$$

\begin{cor} If $X$ is a spherical variety, then $\IH^*_G(X)$ is pure and
$$\IH^n(X)=\bigoplus_{q-p=n} \Tor^{q,H^*(BG)}_p(\IH^*_G(X),\Q)\,.$$
The sum of terms with $q\leq \nu$ coincides with $W_\nu \IH^*(X)$.\end{cor}

This theorem
generalizes a
result for toric varieties,
\cite{We2}.

\section{Koszul duality}

Our paper is motivated by the work of Goresky, 
Kottwitz and MacPherson, \cite{GKM}. One of their
main results states that the 
nonequivariant ``cohomology'' of a $G$-space $X$ can be recovered from 
the
``equivariant cohomology'' (and vice versa) through Koszul duality.
But one
has to be careful, since here by ``equivariant cohomology'' we mean 
a complex $C^*_G(X)$ (an object in a derived category)
computing $H^*_G(X)$. Precisely, there is an equivalence of derived categories
$$D\left(\hbox{$H^*(BG)$-modules}\right)\simeq 
D\left(\hbox{$H_*(G)$-modules}\right)$$
such that $C^*_G(X)$ corresponds to $C^*(X)$.
The argument of \cite{GKM} contains a gap. A proof that the 
cohomology of the Koszul complex $\Omega^*_G(X)\otimes H^*(G)$ is equal 
to 
$H^*(X)$ appeared in $\cite{MW}$\footnote{The action of $H^*(G)$ in
\cite{MW} is not correct.},
while the correct action of $H_*(G)$ has been constructed
in \cite{H}~and~\cite{AM}.
The question arises: can one recover $H^*(X)$
knowing only the cohomology groups $H^*_G(X)$, not the whole complex?
In general one cannot.
The higher differentials can be expressed in terms of Massey products
in $H^*_G(X)$.

\begin{re}\rm A better question seems to be:
Can one find interesting spaces for which knowledge of~$H^*_G(X)$
suffices to recover $H^*(X)$?
This property deserves to be called ``formality'' (or 
maybe $BG$-formality), but unfortunately the word ``formality'' has 
been reserved for something else (namely freeness over~$H^*(BG)$).
See the remarks in the 
introduction to \cite{We2}.\end{re}

In
general there is a 
spectral sequence converging to $H^*(X)$  with
$$E_1=H^*_G(X)\otimes H^*(G)\,,$$
$$d_1=\rm Koszul\; differential\,.$$
This is just the Eilenberg-Moore spectral sequence. The cohomology of 
this complex is a torsion product:
$$E_2={\rm Tor}^{H^*(BG)}(H^*_G(X),\C)\,.$$
Our Theorem \ref{B} says that under suitable purity assumptions the
higher differentials vanish.

\begin{cor} Suppose that $H^*_G(X)$ is pure
then  the cohomology
of $X$ is the cohomology of the Koszul complex $(H^*_G(X)\otimes
H^*(G),d_1)$. \end{cor}

There are
still open questions:
While Koszul duality allows
one to recover the $H_*(G)$-module structure of $H^*(X)$, the
Eilenberg-Moore spectral sequence gives $H^*(X)$ only additively.

\begin{qu}  Suppose that a smooth
complex
 algebraic variety $X$ consists of
finitely many orbits. Is the cohomology of the Koszul complex
$$
(H_G^*(X)\otimes H^*(G); d_1).
$$
isomorphic to $H^*(X)$ as a {\em module over $H_*(G)$}?\end{qu}
By \cite{FT}
this is the case for toric varieties.

Matthias Franz

Section de Math\'ematiques, Universit\'e de Gen\`eve

CP~240, 1211 Gen\`eve 24, Switzerland

e-mail:~\texttt{matthias.franz@math.unige.ch}

\vskip1ex

Andrzej Weber

Instytut Matematyki, Uniwersytet Warszawski

ul. Banacha 2, 02--097 Warszawa, Poland

e-mail:~\texttt{aweber@mimuw.edu.pl}


\begin{thebibliography}{15}
\small
\def\by{ \sc}
\def\book{: \rm}
\def\paper{: \em}
\def\ed{, \rm ed.}
\def\jour{. \rm}
\def\inbook{. \rm}
\def\publ{\rm, }
\def\vol{, \bf}
\def\no{, \rm no.~}
\def\yr#1{ \rm (#1)}
\def\pages{, \rm p.~}
\def\endref{\rm}

\bibitem{AM}\by A.~Alekseev, E.~Meinrenken\paper
Equivariant cohomology and the Maurer-Cartan equation\jour
Preprint math.DG/0406350 
\endref

\bibitem{BBD}\by A. Beilinson, J. Bernstein, P. Deligne\paper Faisceaux
Pervers\jour Ast\'erisque\vol 100\yr{1983} \pages 5--171\endref

\bibitem{BL}\by J.~Bernstein, V.~Lunts\book
Equivariant sheaves and functors, 
Lecture Notes in Mathematics\vol 1578\publ
Springer-Verlag, Berlin, 1994\endref

\bibitem{BDCP}\by E. Bifet, C. De\,Concini, C. Procesi%
\paper Cohomology of Regular Embeddings\jour Adv. Math.\vol 82\no 1\yr{1990}%
\pages 1--34.

\bibitem{Bo}\by A. Borel\paper Sur la cohomologie des espaces fibr\'es et des 
espaces
homog\`enes de groupes de Lie compacts\jour Ann. Math\vol 57\no 1\yr
{1953}\pages 115--207\endref%

\bibitem{Br}\by M. Brion\paper Vari\'et\'es sph\'eriques%
\jour notes, available at\\
{\tt http://www-fourier.ujf-grenoble.fr/\~{ }brion/}

\bibitem{BJ1}\by M. Brion, R. Joshua\paper Vanishing of odd intersection
cohomology II\jour Math. Ann.\vol 321\yr{2001}\pages 399--437\endref  

\bibitem{BJ2}\by M. Brion, R. Joshua\paper Intersection cohomology of
reductive varieties\jour 
arXiv preprint math.AG/0310107\yr{2003} 
\endref  


\bibitem {D1}\by  P.~Deligne\paper
{\em La conjecture de Weil. II}\jour
Inst. Hautes \'Etudes Sci. Publ. Math.\no 52
\yr{1980}\pages137--252\endref


\bibitem{D2}\by P. Deligne\paper Th\'eorie de Hodge III\jour
Publ. Math. I.H.E.S.\vol 44\yr{1974}\pages 5--77\endref

\bibitem{DL}\by J.~Denef, F.~Loeser\paper Weights of exponential
sums, intersection cohomology, and Newton polyhedra\jour Inv. Math.\vol
109\yr{1991}\pages 275--294\endref

\bibitem{EG}\by D.~Edidin, W.~Graham\paper Equivariant
intersection theory\jour Invent. Math.\vol 131\yr{1998)}\no 3\pages
595--634\endref 

\bibitem{EM}\by S. Eilenberg, J.~C.~Moore\paper Homology and fibrations I.
Coalgebras, cotensor product and its derived functors\jour Comment.
Math. Helv.\vol 40\yr{1966}\pages 199--236\endref


\bibitem{FT}\by M.~Franz\paper On the integral cohomology of
smooth toric varieties\jour Preprint math.AT/0308253

\bibitem{GM}\by M. Goresky, R. MacPherson\paper Intersection
homology II\jour Invent. Math.\vol 72\yr{1983}\pages 77--130\endref

\bibitem{GKM}\by M. Goresky, R. Kottwitz, R. MacPherson\paper
Equivariant cohomology, Koszul duality and the localization
theorem\jour Inv. Math.\vol 131\yr {1998}\pages 25--83\endref


\bibitem{Hain}\by R.~M.~Hain\paper The de Rham homotopy theory
of complex algebraic
varieties. I\jour $K$-Theory\vol 1 \yr{1987}\no 3\pages 271--324
\endref

\bibitem{Hain-Zucker}\by R.~M.~Hain, S.~Zucker\paper Unipotent variations
of mixed Hodge
structure\jour Invent. Math.\vol 88\yr{1987}\no 1\pages 83--124
\endref

\bibitem{H}\by J.~Huebschmann\paper
Relative homological algebra, homological perturbations,
equivariant de Rham theory, and Koszul duality\jour 
arXiv preprint math.DG/0401161\yr{2003} 
\endref


\bibitem{MW}\by T.~
Maszczyk, A.~Weber\paper Koszul duality for modules over Lie 
algebras\jour Duke Math. J.\vol 112\yr{2002}\no 3\pages 511--520\endref

\bibitem{MV}\by F.~Morel, V. Voevodsky%
\paper ${\bf A}^1$-homotopy theory of schemes%
\jour Publ. Math. I.H.E.S.\vol 90\yr{2001}\pages45--143.
\endref


\bibitem{Sa1}\by M.~Saito\paper Hodge structure via filtered $\cal
D$-modules\jour Ast\'erisque\vol 130\yr{1985}\pages 342--351
\endref

\bibitem{Sa2}\by M.~Saito\paper Introduction to mixed
Hodge modules\jour
Ast\'erisque\vol 179-180\yr{1989}\pages 145--162.

\bibitem{Sa3}\by M.~Saito\paper Decomposition theorem for proper
K\"ahler morphisms\jour T\^ohoku Math. J.\vol 42 \yr{1990}\pages
127-148.


\bibitem{Sm1}\by L.~Smith\paper On the construction of the
Eilenberg--Moore spectral sequence\jour Bull. Amer. Math. Soc.\vol
75\yr{1969}\pages 873--878\endref
\endref

\bibitem{Sm2}\by L.~Smith\book Lectures on the Eilenberg--Moore
spectral sequence, LNM\vol 134\publ Springer\yr{1970}
\endref

\bibitem{SH}\by J.~Stasheff, S.~Halperin%
\paper Differential algebra in its own rite%
\jour Proc. Adv. Study Inst. Alg. Top. (Aarhus 1970), Vol. III,
  Various Publ. Ser.\vol 13\yr{1970}\pages 567--577
\endref

\bibitem{To}\by B.~Totaro\paper  The Chow ring of a classifying space%
\inbook Algebraic $K$-theory. Proceedings of
an AMS-IMS-SIAM summer research conference, Seattle (WA), USA, July 13--24, 1997.
{\rm Symp. Pure Math. Vol 67}\ed~Raskind, Wayne et al.\publ
Providence, RI: American Mathematical Society. Proc.\pages
249--281\yr{1999}\endref 

\bibitem{We}\by A.~Weber\paper
Formality of equivariant intersection
cohomology of 
algebraic varieties\jour Proc. Amer. Math. Soc.\vol 131\yr{2003}\pages
2633--2638 \endref

\bibitem{We2}\by A.~Weber\paper Weights in the cohomology of toric
varieties\jour 
Central European Journal of Mathematics
\vol 2\yr{2004}\no 3\pages 478--492\endref
\end{thebibliography}
\end{document}